\RequirePackage{lineno}
\documentclass[12pt,a4paper,leqno]{amsart}
\newcounter{minutes}
\divide\time by 60
\newcounter{hours}
\multiply\time by 60
\addtocounter{minutes}{-\time}
 \def\registered{
 {\ooalign{\hfil\raise .00ex\hbox{\scriptsize R}\hfil\crcr\mathhexbox20D}}}

\usepackage{amssymb}
\usepackage{graphicx}
\usepackage{here}
\date{}
\newfont{\cyrilic}{wncyr10 scaled 1000}
\newcommand{\comment}[1]{}

\swapnumbers
\theoremstyle{plain}

\newtheorem{theorem}[equation]{Theorem}

\newtheorem{lemma}[equation]{Lemma}

\newtheorem{mconj}[equation]{The Mori Conjecture}

\newtheorem{corollary}[equation]{Corollary}

\newtheorem{remark}[equation]{Remark}

\newtheorem{subsec}[equation]{}

\numberwithin{equation}{section}          

\font\fFt=eusm10 
\font\fFa=eusm7  
\font\fFp=eusm5  
\def\K{\mathchoice
{\hbox{\,\fFt K}}
{\hbox{\,\fFt K}}
{\hbox{\,\fFa K}}
{\hbox{\,\fFp K}}}

\begin{document}

\thispagestyle{empty}

\vspace{1.0cm}
\title[On Mori's theorem for quasiconformal maps]{\large \bf On Mori's theorem for \\quasiconformal maps in the $n$-space}

\def\thefootnote{}
\footnotetext{
\texttt{\tiny File:~\jobname .tex,
          printed: \number\year-0\number\month-0\number\day,
          \thehours.\ifnum\theminutes<10{0}\fi\theminutes}
}
\makeatletter\def\thefootnote{\@arabic\c@footnote}\makeatother


\author[B. A. Bhayo]{B.A. Bhayo}
\address{Department of Mathematics, University of Turku,
FI-20014 Turku, Finland} \email{barbha@utu.fi}

\author[M. Vuorinen]{M. Vuorinen}
\address{Department of Mathematics, University of Turku,
FI-20014 Turku, Finland} \email{vuorinen@utu.fi}



\maketitle
\bigskip
{\sc In memoriam: M.K. Vamanamurthy, 5 September 1934-- 6 April 2009}
\bigskip

\begin{abstract} R. Fehlmann and M. Vuorinen proved in 1988 that Mori's
constant $M(n,K)$ for $K$-quasiconformal maps of the unit ball in
$\mathbf{R}^n$ onto itself keeping the origin fixed satisfies
$M(n,K) \to 1$ when $K\to 1\,.$ We give here an alternative proof of
this fact, with a quantitative upper bound for the constant in terms
of elementary functions. Our proof is based on a refinement of a
method due to G.D. Anderson and M. K. Vamanamurthy. We also give an
explicit version of the Schwarz lemma for quasiconformal self-maps
of the unit disk. Some experimental results are provided to compare the various bounds
for the Mori constant when $n=2\,.$
\end{abstract}

\bigskip

\section{\sc Introduction}

Distortion theory of quasiconformal and quasiregular mappings in the
Euclidean $n$-space $\mathbf{R}^n$ deals
with estimates for the modulus of continuity and change of distances
under these mappings. Some of the examples are
the H\"older continuity, the quasiconformal counterpart of the Schwarz lemma, and
Mori's theorem. The investigation of these topics started in
the early 1950's for the case $n=2$ and ten years later for the case $n \ge 3 \,.$
Many authors have contributed to the distortion theory, for some historical
remarks see \cite[11.50]{vuobook}.

 As in \cite{fv} we define Mori's constant
$M(n,K)$ in the following way. Let $QC_K,\, K\ge 1,$ stand for the
family of all $K$-quasiconformal maps of the unit ball $
\mathbf{B}^n$ onto itself keeping the origin fixed. Note
that it is a well-known basic fact that an element in the set $QC_K$
can be extended by reflection to a $K$-quasiconformal map of the
whole space $\overline{\mathbf{R}}^n = {\mathbf{R}}^n \cup \{ \infty
\}$ onto itself keeping the point $\infty$ fixed. Then for all $K\ge
1,\, n\ge 2\,,$ there exists a least constant $M(n,K)\ge 1$ such
that
\begin{equation}
 |f(x)-f(y)| \le M(n,K)|x-y|^{\alpha}, \quad \alpha= K^{1/(1-n)}\,,
\end{equation}
for all $ f \in QC_K, x,y\in \mathbf{B}^n$\,(see \cite{fv}).

L. V. Ahlfors \cite{a} proved in 1954 that $M(2,K)\le 12^{K^2}$ and this property was refined by A. Mori
\cite{m} in 1956 to the effect that
 $M(2,K)\leq 16$  and $16$ cannot be replaced by a smaller constant
 independent of $K\,.$ This result can also be found in \cite{a2},
 \cite{fm}, and \cite{lv}.
 On the other hand the trivial observation that $16$
fails to be a sharp constant for $K =1$ led to the following
conjecture, which is still open in 2009.

\begin{mconj} \label{moriconj}   $M(2,K) = 16^{1-1/K}.$
\end{mconj}

 O. Lehto and K.I. Virtanen demonstrated in 1973 \cite[pp. 68]{lv} that
$M(2,K)\geq 16^{1-1/K}$ (this lower bound was not given in the 1965 German
edition of the book). It is natural to expect that for a fixed $n
\ge 2,$ $M(n,K) \to 1$ when $K \to 1$ and this convergence result with an explicit upper bound for $M(n,K)$ was proved by R.
Fehlmann and M. Vuorinen \cite{fv}. A counterpart of this result for the
chordal metric was proved recently by P. H\"ast\"o in \cite{h}.

\begin{theorem}\label{fvu}\cite[Theorem 1.3]{fv} Let $f$ be a $K$-quasiconformal mapping
of $\mathbf{B}^n$ onto $\mathbf{B}^n$, $n\geq 2$, $f(0)=0$. Then
\begin{equation}\label{fvu1}
|f(x)-f(y)| \le M(n,K)|x-y|^{\alpha}
\end{equation}
for all $x,y\in \mathbf{B}^n$ where $\alpha= K^{1/(1-n)}$ and the constant $M(n,K)$
has the following three properties:
\begin{enumerate}
\item $M(n,K)\to 1$ as $K\to 1$, uniformly in $n$\,,
\item $M(n,K)$ remains bounded for fixed $K$ and varying $n$\,,
\item $M(n,K)$ remains bounded for fixed $n$ and varying $K$\,.
\end{enumerate}
\end{theorem}

For $n=2\,,$ the first majorants with the convergence property in
\ref{fvu}(1) were proved only in the mid 1980s and for $n \geq 3$
in \cite{fv}. In \cite{fv} a survey of the various known bounds for
$M(n,K)$ when $n\ge 2$ can be found -- that survey reflects what was
known at the time of publication of \cite{fv}. Some earlier results
on H\"older continuity had been proved in
\cite{g},  \cite{mrv}, \cite{r}, \cite{s}. Step by step the
bound for Mori's constant was reduced during the past twenty years.
As far as we know, the best upper bound known today for $n=2$ is
$M(2,K)\le 46^{1-1/K} $ due to S.-L. Qiu \cite{q} (1997). Refining
the parallel work \cite{fv}, G.D. Anderson and M. K. Vamanamurthy
proved the following theorem in \cite{av}.

\begin{theorem} \label{anva} For $n\geq 2,K\geq 1$,
$$M(n,K) \leq 4\lambda_n^{2(1-\alpha)} \,,   $$
where $\alpha=K^{1/(1-n)}\,$ and $ \lambda_n\in [4, 2e^{n-1}) \,,
\lambda_2=4,$ is the Gr\"otzsch ring constant \cite{an},
\cite[p.89]{vuobook}.
\end{theorem}

The first main result of this paper is Theorem \ref{bvthm} which improves on
Theorem \ref{anva}.

\begin{theorem} \label{bvthm} (1) For $n\geq 2,K\geq 1$, $M(n,K) \le T(n,K)$
\begin{equation} \label{bv1} T(n,K) = \inf \{ h(t) : t \ge 1  \}\,, \quad h(t)
=(3+\lambda_n^{\beta-1} t^{\beta} ) t^{-\alpha}
\lambda_n^{2(1-\alpha)},\;t\ge 1\,,
\end{equation}
where $\alpha=K^{1/(1-n)} =1/\beta,$ and $\lambda_n$ is as in
Theorem \ref{anva}.

(2) There exists a number $K_1>1$ such that for all
$K \in (1,K_1)$ the function $h$ has a minimum at a point $t_1$ with $t_1>1$
and

\begin{equation} \label{bv2}
 T(n,K) \le h(t_1) = \left[ \frac{3^{1-\alpha^2} (\beta-\alpha)^{\alpha^2}}{\alpha^{\alpha^2}} \lambda_n^{\alpha- \alpha^2}
  + \lambda_n^{\beta-1} \left(\frac{(3 \alpha)^\alpha \lambda_n^{\alpha-1}}{(\beta-\alpha)^{\alpha} }   \right) ^{\beta-\alpha} \right]\lambda_n^{2(1-\alpha)} \,.
 \end{equation}
Moreover, for $\beta\in(1,\min\{2,K_1^{1/(n-1)}\})$ we have
\begin{equation}\label{uph}
h(t_1)\leq 3^{1-\alpha^2}2^{5(1-\alpha)}K^5\left(\frac{3}{2}\sqrt[4]{\beta-\alpha}+\exp(\sqrt{\beta^2-1})\right).
\end{equation}
In particular, $h(t_1)\to 1$ when $K\to 1\,.$
\end{theorem}

The last statement shows that Theorem \ref{bvthm}
is better than the result of Anderson and Vamanamurthy,
Theorem  \ref{anva}, at least for values of $K$ close to the critical
value $1$, because the constant of Theorem \ref{anva} satisfies
$4\lambda_n^{2(1-\alpha)} \ge 4. $

The main method of our proof is to replace the argument of Anderson
and Vamanamurthy by a more refined inequality from \cite{vu2} and to
introduce an additional parameter ($t$ in the above theorem) which
will be chosen in an optimal way. The fact that this refined
inequality is essentially sharp for values of $t$ large enough, was
recently proved by V. Heikkala and M. Vuorinen in \cite{hv}. This
gave us a hint that the inequality from \cite{vu2} might lead to an
improvement of the results in \cite{av}. For the case $n=2$ a
numerical comparison of our bound (\ref{bv2}) to Mori's conjectured
bound, to the bound in Theorem \ref{anva} and to the bound in
\cite{fv} is presented in tabular and graphical form at the end of the paper.

We conclude this paper by discussing the Schwarz lemma for plane
quasiconformal self-mappings of the unit disk, formulated in terms
of the hyperbolic metric. The long history of this result is
summarized in \cite[p.152, 11.50]{vuobook}. An up-to-date form of the Schwarz lemma was
given in \cite[Theorem 11.2]{vuobook} and it will be stated for convenient reference
also below as Theorem \ref{my11.2}.
A particular case, formula (\ref{my11.2(2)}), was rediscovered by
D.B.A. Epstein, A. Marden and V. Markovic \cite[Thm 5.1]{emm}.

We use the notations  {\rm ch}, {\rm th}, {\rm arch} and {\rm arth} as in  \cite{vuobook}, to
denote the hyperbolic cosine, tangent and their inverse functions, resp. The
second main result of this paper is an explicit form of the Schwarz lemma for
quasiregular mappings, Theorem \ref{refschlem}. We believe that in this simple form the
result is new and perhaps of independent interest. The constant $c(K)$ below involves the
transcendental function $\varphi_K$ defined in Section 4.

\begin{theorem} \label {refschlem} If $f:\mathbf{B}^2\to \mathbf{R}^2$ is a non-constant
$K$-quasiregular mapping with $f\mathbf{B}^2\subset\mathbf{B}^2$,
and $\rho$ is the hyperbolic metric of $\mathbf{B}^2\,,$ then
$$\rho(f(x),f(y))\leq c(K)\max\{\rho(x,y),\rho(x,y)^{1/K}\}$$
for all $x,y\in \mathbf{B}^2$ where $c(K)=2{\rm arth}(\varphi_K({\rm
th}\frac{1}{2}))\,$ and
$$K\leq u(K-1)+1\leq\log({\rm ch}(K {\rm arch}(e)))\leq c(K)\leq v(K-1)+K$$
with $u={\rm arch}(e){\rm th}({\rm arch}(e))>1.5412$ and $v=\log(2(1+\sqrt{1-1/e^2}))<1.3507$.
In particular, $c(1)=1\,.$
\end{theorem}

{\sc Acknowledgments.} The first author is indebted to the Graduate
School of Mathematical Analysis and its Applications for support.
Both authors wish to acknowledge the kind help of Prof. G.D.
Anderson in the proof of Lemma \ref{anmono}, the valuable help of the
 referee for the improvement of the manuscript, as well as the expert help of Dr. H. Ruskeep\"a\"a in the use of Mathematica$^\registered$ \cite{ru}.

 \vspace{5em}

\section{\sc The main results}

We shall follow here the standard notation and terminology for
$K$-quasiconformal and $K$-quasiregular mappings in the Euclidean $n$-space $
\mathbf{R}^n\,,$ see e.g. \cite{v}, \cite{vuobook}, and we also recall some
basic notation. For the modulus $M( \Gamma)$ of a curve family
$\Gamma$ and its basic properties see \cite{v} and \cite{vuobook}.

Let $D$ and $D^{'}$ be domains in $\overline{\mathbf{R}}^n,K\geq 1$, and let
$f:D\to D^{'}$ be a homeomorphism. Then $f$ is $K$-quasiconformal if
$$M(\Gamma)/K\leq M(f\Gamma)\leq K M(\Gamma)$$
for every curve family $\Gamma$ in $D$ \cite{v}.\

For subsets $E,F, D\subset \overline{\mathbf{R}}^n$ we denote by
$\Delta(E,F;D)$ the family of all curves joining $E$ and $F$ in $D$.
For brevity we write $\Delta(E,F)= \Delta(E,F;{\mathbf{R}}^n)\,.$
A ring is a domain in ${\mathbf{R}}^n$, whose complement consists of
two compact and connected sets. If these sets are $E$ and $F$, then the
ring is denoted by $R(E,F)\,.$ The capacity of a ring $R(E,F)$ is
$${\rm cap}R(E,F)=M(\Delta(E,F)).$$
The complementary components of the Gr\"otzsch ring $R_{G,n}(s)$ are $\overline{\mathbf{B}}^n$
and $[se_1,\infty], s>1$, while those of the Teichm\"uller ring $R_{T,n}(t)$ are $[-e_1,0]$
and $[te_1,\infty],t>0$. The  conformal capacities of $R_{G,n}(s)$ and $R_{T,n}(t)$ are denoted by
$$\left\{\begin{array}{lll}\gamma_n(s)=M(\Delta(\overline{\mathbf{B}}^n,[se_1,\infty]))\\

\tau_n(t)=M(\Delta([-e_1,0],[te_1,\infty]))\end{array}\right.$$
respectively. Here $\gamma_n:(1,\infty)\to(0,\infty)$ and $\tau_n:(0,\infty)\to(0,\infty)$
are decreasing homeomorphisms and they satisfy the fundamental identity

\begin{equation}\label{1}
\gamma_n(s)=2^{n-1}\tau_n(s^2-1),\quad t>1\,,
\end{equation}
see e.g. \cite[5.53]{vuobook}.

For $n\geq 2$ and $K > 0$, the distortion function $\varphi_{K,n}:
[0,1]\to [0,1]$ is a
homeomorphism. It is defined by
\begin{equation}\label{2}
\varphi_{K,n}(t)=\displaystyle\frac{1}{\gamma_n^{-1}(K\gamma_n(1/t))},\quad t\in(0,1),
\end{equation}
and $\varphi_{K,n}(0)=0\,,$  $\varphi_{K,n}(1)=1\,.$  For $n\geq 2, K\geq 1$ and $0\leq r\leq 1$
\begin{equation}\label{3}
\varphi_{K,n}(r)\leq \lambda_n^{1-\alpha}r^\alpha, \quad \alpha=K^{1/(1-n)}\,,
\end{equation}
\begin{equation}\label{3a}
\varphi_{1/K,n}(r)\geq \lambda_n^{1-\beta}r^\beta,
\quad \beta=K^{1/(n-1)} \,,
\end{equation}
by \cite[Theorem 7.47]{vuobook} and where $\lambda_n\ge 4$ is as in
Theorem \ref{anva}.

\begin{lemma}\label{3b} Suppose that $f : \mathbf{B}^n\to\mathbf{B}^n$ is a $K$-quasiconformal mapping with
 $f\mathbf{B}^n=\mathbf{B}^n$, $f(0)=0,$ and let $h : \overline{\mathbf{R}}^n \to \overline{\mathbf{R}}^n$
 be the inversion $h(x) = x/|x|^2\,, h(\infty)=0, h(0)=\infty,$ and define
 $g:\overline{\mathbf{R}}^n\to\overline{\mathbf{R}}^n$
 by $g(x)=f(x)$ for $x\in \mathbf{B}^n, g(x)=h(f(h(x)))$ for $x\in\mathbf{R}^n\setminus
 \overline{\mathbf{B}}^n$ and $g(x)=\lim_{z\to x}f(z)$ for $x\in\partial\mathbf{B}^n,g(\infty)=\infty$.
 Then $g$ is a $K$-quasiconformal mapping, and we have for $x\in\mathbf{B}^n$
\begin{equation}\label{3c}
\varphi_{1/K,n}(|x|)\leq |f(x)|\leq\varphi_{K,n}(|x|).
\end{equation}
For $x\in\mathbf{R}^n\setminus\overline{\mathbf{B}}^n$
\begin{equation}\label{3d}
1/\varphi_{K,n}(1/|x|)\leq |g(x)|\leq1/\varphi_{1/K,n}(1/|x|).
\end{equation}
\end{lemma}
\begin{proof} It is well-known that the above definition defines $g$ as a $K$-quasiconformal
homeomorphism. The formula (\ref{3c}) is well-known (see \cite[Theorem 4.2]{avv}) and (\ref{3d}) follows easily.
\end{proof}

\begin{lemma}\label{4}\cite[Lemma 7.35]{vuobook}
Let $R=R(E,F)$ be a ring in $\overline{\mathbf{R}}^n$ and let $a,b\in E,c,d\in F$ be distinct points. Then
$$\text{cap}R=M(\Delta(E,F))\geq\tau_n\left(\frac{|a-c||b-d|}{|a-b||c-d|}\right).$$
Equality holds if $b=t_1e_1,a=t_2e_1,c=t_3e_1,d=t_4e_1$ and $t_1<t_2<t_3<t_4$.
\end{lemma}

We consider Teichm\"uller's extremal problem, which will be used to
provide a key estimate in what follows. For
$x\in\mathbf{R}^n\setminus\{0,e_1\}, n\geq 2$, define
$$p_n(x)=\inf_{E,F}M(\Delta(E,F))$$
where the infimum is taken over all the pairs of continua $E$ and
$F$ in $\overline{\mathbf{R}}^n$ with $0,e_1\in E$ and $x,\infty \in
F$. Note that Lemma \ref{4} gives the lower bound for $p_n(x)$
in Lemma \ref{5}.

\begin{lemma}\label{5}\cite[Theorem 3.20]{vu2}
For $z\in\mathbf{R}^n, |z|>1$, the following inequalities hold:
$$\tau_n(|z|)=p_n(-|z|e_1)\leq p_n(z)\leq p_n(|z|e_1)=\tau_n(|z|-1)$$
where $p_n(z)$ is the Teichm\"uller function. Furthermore, for
$z\in\mathbf{R}^n\setminus[0,e_1]$, there exists a circular arc $E$
with $0,e_1\in E$ and a ray $F$ with $z,\infty\in F$
such that

\begin{equation}\label{6}
p_n(z)\leq \tau_n\left(\frac{|z|+|z-e_1|-1}{2}\right)=M(\Delta(E,F))
\leq \tau_n(|z|-1)
\end{equation}
with equality in the first inequality both for $z=-se_1, s>0$, and for $z=se_1,s>1\,.$
\end{lemma}

\begin{subsec}{\bf Notation.}\label{6b} {\rm For $t>0, x,y \in \mathbf{B}^n\,,$ we write
$$ D(t,x,y)=|x+t\frac{y}{|y|}| \,\,\, \mathrm{if \,\,} y \neq 0,\quad  D(t,x,0)=|x+ e_1|\,.$$
By the triangle inequality we have}
\begin{equation}\label{6c}
t- |x|\le D(t,x,y) \le t+ |x|\,.
\end{equation}
\end{subsec}

\begin{theorem}\label{7a} For $n\geq 2,K\geq 1$, let $f:\overline{\mathbf{R}}^n\to\overline{\mathbf{R}}^n$
be a $K$-quasiconformal mapping, with
$f\mathbf{B}^n=\mathbf{B}^n$, $f(0)=0$ and $f(\infty)=\infty$.
Then for $t\ge 1\,,$ $x,y\in \mathbf{B}^n\setminus \{0\}\,,$ we have
\begin{eqnarray*}
|f(x)-f(y)|&\leq &(3+\varphi_{1/K,n}(1/t)^{-1})
\varphi_{K,n}^2\left(\left(\frac{2|x-y|}{s_1+|x-y|}\right)^{1/2}\right)\\
          &\leq & (3+\lambda_{n}^{(\beta-1)}t^\beta)\lambda_n^{2(1-\alpha)}\left(\frac{2|x-y|}{s_1+|x-y|}\right)^\alpha\,,\;
          \alpha=K^{1/(1-n)}=1/\beta,
\end{eqnarray*}
where $s_1=\displaystyle\max\{a,b \}, a=t+|x|+ D(t,y,x),
b=t+|y|+D(t,x,y)$. 
\end{theorem}

\begin{proof}
Let $\Gamma$ be the family $\Delta(E,F)$ and let
$E$ and $F$ be
 connected sets as in Lemma \ref{5} with
$x,y\in E,z,\infty \in F$, where $z=-tx/|x|$ and $\Gamma^{'}=f(\Gamma)$.
By Lemma $\ref{4}$ and ($\ref{6}$), we have
\begin{eqnarray*}
& {} &\tau_n\left(\frac{|f(z)-f(x)|}{|f(x)-f(y)|}\right)\leq  M(\Gamma^{'})\leq KM(\Gamma)\\
                                                  &\leq & K\tau_n(u-1)\,, \quad
u=\displaystyle\frac{|x-z|+|z-y|+|x-y|}{2|x-y|}\,.
\end{eqnarray*}
The basic identity ($\ref{1})$ yields
\begin{equation}\label{8}
 \gamma_n\left(\left(\displaystyle\frac{|f(z)-f(y)|+|f(x)-f(y)|}{|f(x)-f(y)|}\right)^{1/2}\right)\leq
K\gamma_n\left(\left(u\right)^{1/2}\right)
\end{equation}
$$
=K\gamma_n\left(\left(\displaystyle\frac{t+|x|+D(t,y,x)+|x-y|}{2|x-y|}\right)^{1/2}\right).
$$
Applying $\gamma_n^{-1}$ to (\ref{8}) we have

\begin{figure}
\includegraphics[width=8cm]{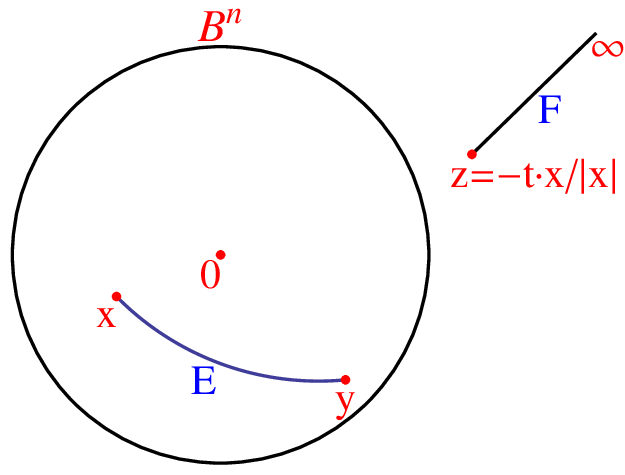}
\caption{ Geometrical meaning of the proof of Theorem \ref{7a}.}
\end{figure}

$$\displaystyle\frac{|f(z)-f(y)|+|f(x)-f(y)|}{|f(x)-f(y)|}\geq
\displaystyle\left(\gamma_n^{-1}\left(K\gamma_n\left(\left
(\frac{a+|x-y|}{2|x-y|}\right)^{1/2}\right)\right)\right)^2=v.$$
Because $f\mathbf{B}^n=\mathbf{B}^n$, by ($\ref{3c}$) and ($\ref{3a}$) we know that
$$|f(z)-f(y)|+|f(x)-f(y)|\leq 3+\varphi_{1/K,n}(1/t)^{-1}\leq 3+\lambda_{n}^{(\beta-1)}t^\beta,$$
\begin{equation}
\frac{|f(x)-f(y)|}{3+ \varphi_{1/K,n}(1/t)^{-1}} \leq  \frac{|f(x)-f(y)|}{|f(z)-f(y)|+|f(x)-f(y)|}
                     \leq 1/v,
\end{equation}  
also
\begin{eqnarray*}
|f(x)-f(y)|&\leq & (3+\varphi_{1/K,n}(1/t)^{-1}) \varphi_{K,n}^2\left(\left(\frac{2|x-y|}{a+|x-y|}\right)^{1/2}\right)\\
            &\leq & (3+\lambda_{n}^{(\beta-1)}t^\beta)\lambda_n^{2(1-\alpha)}\left(\frac{2|x-y|}{a+|x-y|}\right)^\alpha
\end{eqnarray*}
by inequalities (\ref{2}) and (\ref{3}). Exchanging the roles of $x$
and $y$ we see that
\begin{eqnarray*}
|f(x)-f(y)|&\leq &(3+\varphi_{1/K,n}(1/t)^{-1})
\varphi_{K,n}^2\left(\left(\frac{2|x-y|}{s_1+|x-y|}\right)^{1/2}\right)\\
          &\leq & (3+\lambda_{n}^{(\beta-1)}t^\beta)\lambda_n^{2(1-\alpha)}\left(\frac{2|x-y|}
          {s_1+|x-y|}\right)^\alpha.
\end{eqnarray*}
\end{proof}

Setting $t=1$, we get the following corollary.
\begin{corollary}\label{7} For $n\geq 2,K\geq 1$, let $f:\overline{\mathbf{R}}^n\to\overline{\mathbf{R}}^n$
be a $K$-quasiconformal mapping, with $f\mathbf{B}^n = \mathbf{B}^n$,
$f(0)=0$ and $f(\infty)=\infty$. Then for all $x,y\in \mathbf{B}^n\setminus \{0\}\,,$
$$|f(x)-f(y)|\leq 4\lambda_n^{2(1-\alpha)}\left(\frac{2|x-y|}{s+|x-y|}\right)^\alpha\,,$$
where $\alpha=K^{1/(1-n)}$ and $s=\displaystyle\max\{a,b\},
a=1+|x|+D(1,y,x), b=1+|y|+D(1,x,y)\,.$
\end{corollary}

\begin{proof} The proof is similar to the above proof except that here we consider the particular case $t=1$.
Because $f\mathbf{B}^n=\mathbf{B}^n$, we know that $|f(z)-f(y)|+|f(x)-f(y)|\leq 4$,
\begin{eqnarray*}
\frac{|f(x)-f(y)|}{4}&\leq & \frac{|f(x)-f(y)|}{|f(z)-f(y)|+|f(x)-f(y)|}\\
                     &\leq & \displaystyle\frac{1}{\left(\gamma_n^{-1}\left(K\gamma_n\left(\left
(\displaystyle\frac{a+|x-y|}{2|x-y|}\right)^{1/2}\right)\right)\right)^2},
\end{eqnarray*}
or
\begin{eqnarray*}
|f(x)-f(y)|&\leq & 4 \varphi_{K,n}^2\left(\left(\frac{2|x-y|}{a+|x-y|}\right)^{1/2}\right)\\
            &\leq & 4\lambda_n^{2(1-\alpha)}\left(\frac{2|x-y|}{a+|x-y|}\right)^\alpha
\end{eqnarray*}
by inequalities (\ref{2}) and (\ref{3}). Exchanging the roles of $x$
and $y$ we get
$$|f(x)-f(y)| \leq  4\lambda_n^{2(1-\alpha)}\left(\frac{2|x-y|}{\max\{a,b\}+|x-y|}\right)^\alpha \,.$$
\end{proof}

\begin{corollary}\label{9a} For $n\geq 2, K\geq 1, t\ge 1$, let $f$ be as in Theorem \ref{7a}. Then
\begin{equation}\label{11a}
|f(x)-f(y)| \leq  (3+\lambda_n^{(\beta-1)}t^\beta)\lambda_n^{2(1-\alpha)}
\left(\frac{2|x-y|}{2t+||x|-|y||+|x-y|}\right)^\alpha,
\end{equation}
for all $x,y\in \mathbf{B}^n\,,$
\begin{equation}\label{11aa}
|f(x)-f(y)| \leq (3+\lambda_n^{\beta-1}t^\beta)\lambda_n^{2(1-\alpha)}
\left(\frac{|x-y|}{\max\{t+|x|,t+|y|\}}\right)^\alpha,
\end{equation}
for all $x,y\in \mathbf{B}^n\,,$ and

\begin{equation}\label{12a}
|f(x)-f(y)| \leq (3+\lambda_n^{(\beta-1)}t^\beta)\lambda_n^{2(1-\alpha)}
\left(\frac{|x-y|}{t+|x|+(|x-y|)/2}\right)^\alpha,
\end{equation}
if $D(t,y,x)>t+|x|, x,y\in \mathbf{B}^n$.
\end{corollary}

\begin{proof} 
Inequality $(\ref{11a})$ follows because by (\ref{6b})
$D(t,y,x)>t-|y|$ and $D(t,x,y)>t-|x|$ for  $x,y \in \mathbf{B}^n$,
and hence, in the notation of Theorem \ref{7a},
$$s_1\geq\max\{2t+|x|-|y|,2t+|y|-|x|\}=2t+||x|-|y||\,.$$
It is also clear that $D(t,y,x)\geq t+|x|-|x-y|$, and this      
implies that
$$s_1\geq \max \{2(t+|x|)-|x-y|,2(t+|y|)-|x-y|\}=
2 \max\{t+|x|,t+|y|\}-|x-y|$$ and hence the inequality $(\ref{11aa})$
follows. In the case of $(\ref{12a})$ we have $D(t,y,x)>t+|x|$ and see that, in the notation of
Corollary \ref{7}, $s>2(t+|x|)$ and $(\ref{12a})$ holds.
\end{proof}

\begin{corollary}\label{13} For $n\geq 2, K\geq 1$, let $f$ be as in Theorem \ref{7a}. Then
\begin{equation}\label{14}
|f(x)-f(y)| \leq  4 \lambda_n^{2(1-\alpha)}
\left(\frac{2|x-y|}{2+||x|-|y||+|x-y|}\right)^\alpha,
\end{equation}
for all $x,y\in \mathbf{B}^n \setminus \{0 \}\,.$
\end{corollary}

\begin{remark}{\rm (1) In several of the above results we have supposed that
$x,y\in \mathbf{B}^n \setminus \{0 \}\,.$ If one of the points $x,y$ were equal to
$0\,,$ then we would have a better result from the Schwarz lemma
estimate (\ref{wang}).

(2) Corollary \ref{13} is an improvement of the Anderson-Vamanamurthy theorem \ref{anva}\,. }
\end{remark}

\section{\sc Comparison with earlier bounds}

\begin{subsec} {\bf Proof of Theorem \ref{bvthm}.} {\rm (1) The inequality (\ref{bv1})
follows easily from the inequality (\ref{11aa}).

(2) We see that the function $h$ has a local minimum at $t_1= (3
\alpha)^{\alpha} \lambda_n^{\alpha-1} (\beta-\alpha)^{-\alpha} \,.$ If $t_1
\ge 1\,,$ then the inequality (\ref{11aa}) yields the desired
conclusion. The upper bound for $T(n,K)$ follows by substituting the
argument $t_1$ in the expression of $h\,.$ 

We next show that the value $K_1=4/3$ will do. Fix $K\in(1,K_1)\,.$
Then $\alpha= K^{1/(1-n)}\ge 3/4$ and $\alpha/(1- \alpha^2)>1$.

Because $\lambda_n^{ \alpha-1}\ge 2^{1/K -1} K^{-1}$ by \cite[Lemma 7.50(1)]{vuobook}, with $d=(6/K)^{1/K}/2K$
we have
$$
t_1=(3\alpha)^\alpha \lambda_n^{\alpha-1}(\beta-\alpha)^{-\alpha} \ge (3/K)^{1/K} 2^{1/K-1} K^{-1}
\left(\frac{\alpha}{1-\alpha^2} \right)^{\alpha}
$$
$$
=d\left(\frac{\alpha}{1-\alpha^2} \right)^{\alpha}\ge d \left(\frac{\alpha}{1-\alpha^2} \right)^{3/4}
$$
$$
  = \left(2 r(K) \frac{\alpha}{1-\alpha^2} \right)^{3/4}\,;
 r(K)= d^{4/3}/2 \,.
$$
It suffices to observe that $t_1>1$ certainly holds if $2r(K)(\frac{\alpha}{1-\alpha^2})>1$ which holds
for $ \alpha > 1/(r(4/3)+  \sqrt{1+r(4/3)^2})= 0.53...\,,$ in particular, $t_1>1$ holds in the present case
$\alpha>3/4\,.$

For the proof of (\ref{uph}) we give the following inequalities
\begin{equation}\label{upha}
\lambda_n^{\alpha- \alpha^2}\leq 2^{\alpha(1- \alpha)}K^\alpha
\leq 2^{1- \alpha}K^\alpha, \quad K\geq 1 \,,
\end{equation}
\begin{equation}\label{uphb}
\lambda_n^{\beta- \alpha}= \lambda_n^{\beta+1-1- \alpha}= \lambda_n^{\beta(1-\alpha)+1- \alpha}
=\lambda_n^{(\beta+1)(1-\alpha)}\leq(2^{1-\alpha}K)^3,\quad \beta\in(1,2) \,,
\end{equation}
see \cite[Lemma 7.50(1)]{vuobook}. The formula (\ref{bv2}) for
$h(t_1)$ has two terms. We estimate separately each term as follows
\begin{eqnarray*}
\frac{3^{1-\alpha^2} (\beta-\alpha)^{\alpha^2}}{\alpha^{\alpha^2}} \lambda_n^{\alpha- \alpha^2}\lambda_n^{2(1-\alpha)}
&\leq& \frac{3^{(1-\alpha)(1+\alpha)}2^{\alpha(1-\alpha)}2^{2(1-\alpha)}K^2 (\beta-\alpha)^{\alpha^2}}{\alpha^{\alpha^2}} K^\alpha\\
&\leq& \frac{(9\cdot 2\cdot 4)^{1-\alpha} K^2(\beta-\alpha)^{\alpha^2}}{\alpha^{\alpha^2}} K^\alpha\\
&=& 72^{1-\alpha}(\beta-\alpha)^{\alpha^2}K^2K^\alpha\exp(-\alpha^2\log\alpha)\\
&\leq & 72^{1-\alpha}(\beta-\alpha)^{\alpha^2}K^2K^{\alpha}\exp(-\alpha\log\alpha)\\
&=& 72^{1-\alpha}(\beta-\alpha)^{\alpha^2}K^2\exp((\log K-\log\alpha)\alpha)\\
&= & 72^{1-\alpha}(\beta-\alpha)^{\alpha^2}K^2\exp\left(\left(1+\frac{1}{n-1}\log K\right)\alpha\right)\\
&=& 72^{1-\alpha}(\beta-\alpha)^{\alpha^2}K^2\exp\left(\frac{n}{n-1}\alpha\log K\right)\\
&\leq& 72^{1-\alpha}(\beta-\alpha)^{\alpha^2}K^2\exp(2\log K)\\
&= & 72^{1-\alpha}(\beta-\alpha)^{\alpha^2}K^4
\end{eqnarray*}
by inequality (\ref{upha}),

\begin{eqnarray*}
\lambda_n^{2(1-\alpha)}\lambda_n^{\beta-1} \left(\frac{(3 \alpha)^\alpha \lambda_n^{\alpha-1}}{(\beta-\alpha)^{\alpha} }\right) ^{\beta-\alpha}
&=&\lambda_n^{2(1-\alpha)}\lambda_n^{\beta-1}\left((3 \alpha)^\alpha \lambda_n^{\alpha-1}\right)^{\beta-\alpha}(\beta-\alpha)^{-\alpha(\beta-\alpha)}\\
&\leq &(2^{1-\alpha}K)^2\lambda_n^{\beta-\alpha}\left((3 \alpha)^\alpha \lambda_n^{\alpha-1}\right)^{\beta-\alpha}
\left(\frac{\beta^2-1}{\beta}\right)^{-\alpha((\beta^2-1)/\beta)}\\
&\leq &(2^{1-\alpha}K)^2\left(3 ^\alpha \lambda_n\right)^{\beta-\alpha}\beta^{\alpha^2}(\beta^2-1)^{-\alpha^2(\beta^2-1)}\\
&\leq &(2^{1-\alpha}K)^2  3^ {\alpha(\beta-\alpha)} \lambda_n^{(\beta+1)(1-\alpha)}\exp\left(\frac{2\alpha^2}{e}\sqrt{\beta^2-1}\right)\\
&\leq &3^{1-\alpha^2}(2^{1-\alpha}K)^2 (2^{1-\alpha}K)^{(\beta+1)}\exp\left(\frac{2\alpha^2}{e}\sqrt{\beta^2-1}\right)\\
&\leq &3^{1-\alpha^2}(2^{1-\alpha}K)^5\exp(\sqrt{\beta^2-1}),
\end{eqnarray*}
here we assume that $\beta\in(1,2)$ which implies that
$\alpha\in(1/2,1)$. Also the inequalities
$(K-1)^{-(K-1)}\leq\exp((2/e)\sqrt{K-1})$ and (\ref{uphb}) were
used, and we get
\begin{equation}\label{c}
h(t_1)\leq\left[ 72^{1-\alpha}(\beta-\alpha)^{\alpha^2}K^4+3^{\beta-\alpha}(2^{1-\alpha}K)^5\exp(\sqrt{\beta^2-1})\right].
\end{equation}

Because $(\beta-\alpha)\in(0,\frac{3}{2})$ this implies that
$\frac{2}{3}(\beta-\alpha)\in(0,1)$ and $\alpha^2\in(\frac{1}{4},1)$
and further $(\frac{2}{3}(\beta-\alpha))^ {\alpha^2}\leq
(\frac{2}{3}(\beta-\alpha))^{1/4}$, and finally
$$(\beta-\alpha)^
{\alpha^2}\leq(2/3)^{-\alpha^2}(\frac{2}{3}(\beta-\alpha))^
{1/4}\leq (3/2)^{3/4}\sqrt[4]{\beta-\alpha}$$
$$=(3/2)^{3/4}\sqrt[4]{\beta-\alpha}<(3/2)\sqrt[4]{\beta-\alpha}\,.$$
Next we prove that
\begin{equation} \label{d}
72^{1-\alpha}\leq 3^{1-\alpha^2}2^{5(1-\alpha)}K\,.
\end{equation}
This inequality is equivalent to
$$2^{2(\alpha-1)}3^{(1-\alpha)^2}\leq K\Longleftrightarrow -(1-\alpha)\log 4+(1-\alpha)^2\log 3\leq \log K\,.$$
This last inequality holds because the left hand side is negative.
Now from (\ref{c}) and (\ref{d}) we get the desired inequality (\ref{uph}). } $\square$
\end{subsec}

\begin{subsec}{\bf Graphical and numerical comparison of various bounds.}
{\rm The above bounds involve the Gr\"otzsch ring constant $\lambda_n,$ which is known only for $n=2, \lambda_2=4.$
Therefore only for $n=2$ we can compute the values of the bounds.
Solving numerically the equation ${4\cdot16}^{1-1/K}=h(t_1) $
for $K$ we obtain $K=1.3089\,.$ We give numerical and graphical
comparison of the various bounds for the Mori constant. }

{\rm Tabulation of the various upper bounds for Mori's constant when
$n=2$ and $\lambda_2=4$ as a function of $K$: (a) Mori's conjectured
bound $16^{1-1/K}$, (b) the Anderson-Vamanamurthy bound $4 \cdot
16^{1-1/K}$, (c) the bound from (\ref{bv2}). For $K\in (1, 1.3089)$
the upper bound in (\ref{bv2}) is better than the
Anderson-Vamanamurthy bound. Note that the upper bound $T(n,K)\le h(t_1)$ in
(\ref{bv2}) is proved only for $K \in (1,K_1), K_1= 4/3\,.$
We do not know whether it holds for larger values of $K$ but just comparing
the values of $h(t_1)$ and the bound of Fehlmann and Vuorinen for
 $K >1.5946$ we see that $h(t_1)$ is the smaller one of these two.
Numerical values of the \cite{fv} bound given in the table were
computed with the help of the algorithm for $\varphi_{K,2}(r)$
attached with {\cite[p. 92, 439]{avvb}}}.

\begin{figure}
\includegraphics[width=12cm]{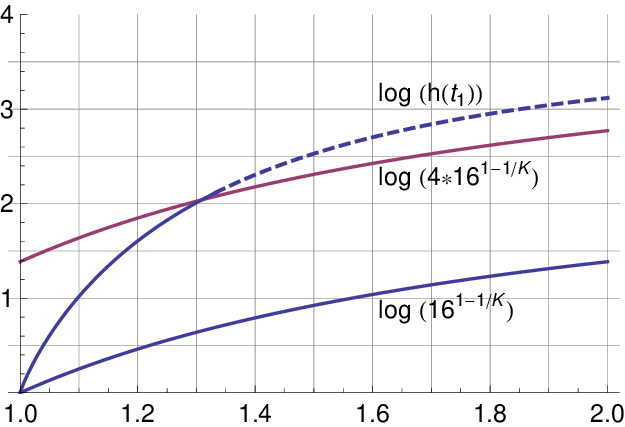}
\caption{Graphical illustration of the various upper bounds for Mori's constant when $n=2$
and $\lambda_2=4$ as a function of $K$:
(a) Mori's conjectured bound $16^{1-1/K}$, (b) the Anderson-Vamanamurthy bound $4 \cdot 16^{1-1/K}$,
(c) the bound from (\ref{bv2}),  valid for $K\in(1,K_1), K_1=4/3$. For $K\in (1, 1.3089)$ the upper bound in (\ref{bv2}) is better
than the Anderson-Vamanamurthy bound.}
\end{figure}

\begin{displaymath}
\begin{array}{|c|c|c|c|c|}
\hline
K&\log({16}^{1-1/K})& \log({4\cdot16}^{1-1/K})&\log(FV)&\displaystyle\log(h(t_1))\\
\hline

 1.1 & 0.2521 & 1.6384 & 0.7051 & 1.0188 \\

 1.2 & 0.4621 & 1.8484 & 1.2485 & 1.6058 \\

 1.3 & 0.6398 & 2.0261 & 1.7046 & 2.0107 \\

 1.4 & 0.7922 & 2.1785 & 2.0913 & 2.3061 \\

 1.5 & 0.9242 & 2.3105 & 2.4221 & 2.5296 \\

 1.6 & 1.0397 & 2.4260 & 2.7094 & 2.7031 \\

 1.7 & 1.1417 & 2.5280 & 2.9633 & 2.8409 \\

 1.8 & 1.2323 & 2.6186 & 3.1921 & 2.9521 \\

 1.9 & 1.3133 & 2.6996 & 3.4020 & 3.0433 \\

 2.0 & 1.3863 & 2.7726 & 3.5979 & 3.1192\\
\hline
\end{array}
\end{displaymath}

{\rm
Note that according to Theorem \ref{bvthm} the inequality (\ref{bv2}) involving
 $h(t_1)$
holds for $\,K\in(1,K_1)$ where the number $K_1>1$ may be smaller than $2.$


For graphing and tabulation purposes we use the logarithmic scale.
Note that the upper bound for $M(2,K)$ given in \cite[Theorem 2.29]{fv}
also has the desirable property that it converges to $1$
when $K\to 1\,,$ see Figure 3.
}

\begin{figure}
\includegraphics[width=12cm]{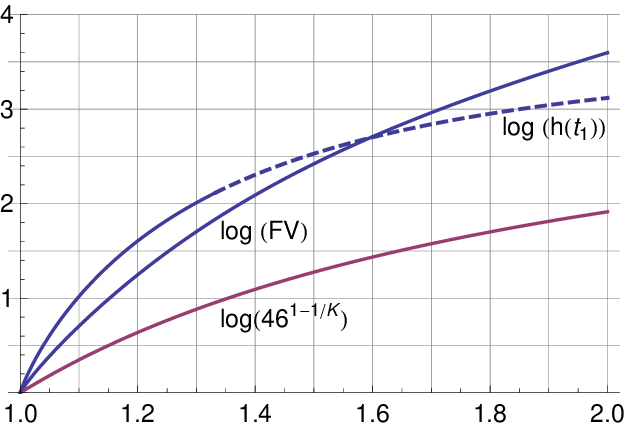}
\caption{Graphical comparison of various bounds
 when $n=2$ and $\lambda_2=4\,,$ as a function of $K$:
(a) the bound from (\ref{bv2}), valid for $K\in(1,K_1), K_1=4/3$,
  (b) the Fehlmann and Vuorinen bound \cite{fv}
$$M(2,K) \le \left(1+\varphi_{K,2}\left(\frac{K^2-1}{K^2+1}\right)\right)2^{2K-3/K}
\frac{(K^2+1)^{(K+1/K)/2}}{(K^2-1)^{(K-1/K)/2}}\,$$
(c) Qiu's bound $46^{1-1/K}$ \cite{q}.}
\end{figure}

\end{subsec}

\begin{subsec}{\bf Comparison of estimates for the H\"older quotient.}
{\rm For a $K$-quasiconformal mapping
$f:\mathbf{B}^n\to f\mathbf{B}^n=\mathbf{B}^n\, , $
we call the expression
$$HQ(f)=\sup\{|f(x)-f(y)|/|x-y|^\alpha:\; x,y\in\mathbf{B}^n, f(0)=0\,\;x\neq y \},$$
the H\"older coefficient of $f$. Clearly $HQ(f)\leq M(n,K)$. Theorem
\ref{7a} yields, after dividing the both sides of the inequality in
\ref{7a} by $|x-y|^{\alpha}\,,$ the upper bound $HQ(f)\le HQ(K)$ for
the H\"older quotient with
\begin{equation}\label{hq}
HQ(K)= \sup\{\inf \{U(t,x,y): \;t\geq 1\}: \; x,y\in\mathbf{B}^n\} \, ,
\end{equation}
$$U(t,x,y)=(3+\varphi_{1/K,n}(1/t)^{-1})
\varphi_{K,n}^2\left(\left(\frac{2|x-y|}{s_1+|x-y|}\right)^{1/2}\right)\frac{1}{|x-y|^\alpha}\,.$$
For $n=2$ we compare $HQ(K)$ to several other bounds (a) Mori's
conjectured bound, (b) the FV bound, (c) the AV bound and give the
results as a table and Figure 3. Because the supremum and infimum in
(\ref{hq}) cannot be explicitly found we use numerical methods that
come with Mathematica software. For the numerical tests we used for
the supremum a sample of $100,000$ random points of the unit disk.}

\bigskip
\begin{figure}
\includegraphics[width=12cm]{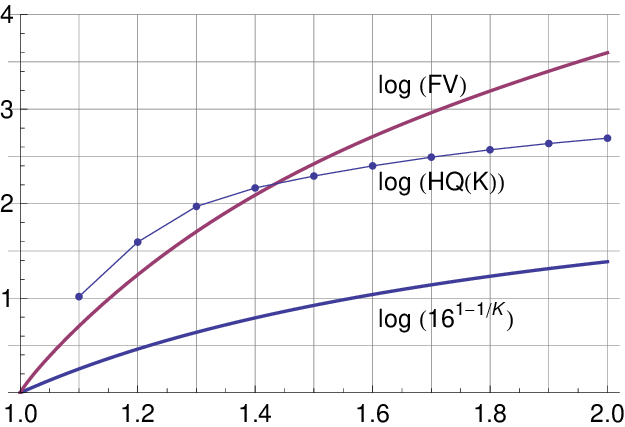}
\caption{Graphical comparison of various bounds
 when $n=2$ and $\lambda_2=4\,,$ as a function of $K$:
(a) the bound from (\ref{hq}),  (b) the Fehlmann and Vuorinen bound \cite{fv}
$$M(2,K) \le \left(1+\varphi_{K,2}\left(\frac{K^2-1}{K^2+1}\right)\right)2^{2K-3/K}
\frac{(K^2+1)^{(K+1/K)/2}}{(K^2-1)^{(K-1/K)/2}},$$ (c) the bound of the Mori conjecture.
The bound  (\ref{hq}) is based on a simulation with  $100,000$ random pairs of points.}
\end{figure}

\begin{displaymath}
\begin{array}{|c|c|c|c|c|}
\hline
K&\log({16}^{1-1/K})& \log({4\cdot16}^{1-1/K})&\log(FV)&\log(HQ(K))\\
\hline

 1.1 & 0.2521 & 1.6384 & 0.7051 & 1.0171 \\

 1.2 & 0.4621 & 1.8484 & 1.2485 & 1.5940 \\

 1.3 & 0.6398 & 2.0261 & 1.7046 & 1.9712\\

 1.4 & 0.7922 & 2.1785 & 2.0913 & 2.1668 \\

 1.5 & 0.9242 & 2.3105 & 2.4221 & 2.2928 \\

 1.6 & 1.0397 & 2.4260 & 2.7094 & 2.4003 \\

 1.7 & 1.1417 & 2.5280 & 2.9633 & 2.4922 \\

 1.8 & 1.2323 & 2.6186 & 3.1921 & 2.5706 \\

 1.9 & 1.3133 & 2.6996 & 3.4020 & 2.6371 \\

 2.0 & 1.3863 & 2.7726 & 3.5979 & 2.6934 \\

\hline
\end{array}
\end{displaymath}

\end{subsec}
\section{An explicit form of the Schwarz lemma}

Recall that the hyperbolic metric $\rho(x,y), x,y \in \mathbf{B}^n\,,$
of the unit ball is given by (cf. \cite{kl},
\cite{vuobook})
\begin{equation} \label{tanhid}
{\rm th}^2\frac{\rho(x,y)}2=\frac{|x-y|^2}{|x-y|^2+t^2}\, ,\quad t^2=(1-|x|^2)(1-|y|^2) \,.
\end{equation}

Next, we consider a decreasing homeomorphism
$\mu:(0,1)\longrightarrow(0,\infty)$ defined by
\begin{equation}
\label{dec_hom} \mu(r)=\frac\pi 2\,\frac{{\K}(r')}{{\K}(r)}, \quad {\K}(r)=\int_0^1\frac{dx}{\sqrt{(1-x^2)(1-r^2x^2)}}\, ,
\end{equation}
where ${\K}(r)$ is Legendre's complete elliptic integral of the
first kind and $ r'=\sqrt{1-r^2},$ for all $r\in(0,1)$.

The Hersch-Pfluger distortion function is an increasing
homeomorphism $\varphi_K:(0,1)\longrightarrow(0,1)$ defined by setting
\begin{equation} \label{phidef}
\varphi_K(r)=\mu^{-1}(\mu(r)/K)\,,\, r\in(0,1),\,\, K>0.
\end{equation}
Note that with the notation of Section
2, $\gamma_2(1/r) = 2 \pi/\mu(r)$ and $\varphi_K(r) =
\varphi_{K,2}(r)\,$ for $r \in (0,1)\,.$

\bigskip

\begin{theorem} \label{my11.2} {\rm \cite[11.2]{vuobook}}{ Let $f :\mathbf{B}^n \to \mathbf{R}^n$ be a
nonconstant  $K$-quasiregular mapping with $ f \mathbf{B}^n \subset
\mathbf{B}^n $, $n\geq 2$, and let $\alpha= K^{1/(1-n)}\,.$ Then }
\begin{equation} \label{my11.2(1)}
{\rm th }\frac{\rho(f(x),f(y))}{2} \le \varphi_{K,n}({\rm th }\frac{\rho(x,y)}{2}) \le
\lambda_n^{1- \alpha} \left({\rm th }\frac{\rho(x,y)}{2}\right)^{\alpha}
\,,
\end{equation}
\begin{equation} \label{my11.2(2)}
\rho(f(x),f(y)) \le K(\rho(x,y) + \log 4)\,,
\end{equation}
{ for all $x,y\in \mathbf{B}^n\,,$ where $\lambda_n$ is the same
constant as in {\rm (\ref{anva})}. If $f(0) = 0\,,$ then}
\begin{equation} \label{wang}
|f(x)| \le  \lambda_n^{1- \alpha} |x|^{\alpha}\,,
\end{equation}
{for all $x\in \mathbf{B}^n\,.$}
\end{theorem}

\bigskip

In the case of quasiconformal mappings with $n=2$ formulas (\ref{my11.2(1)}) and (\ref{wang})
also occur in  \cite[p. 65]{lv} and formula (\ref{my11.2(2)}) was rediscovered in \cite[Theorem 5.1]{emm}.
Comparing Theorem \ref{my11.2} to Theorem \ref{refschlem} we see
that for $n=2$ the expression $K(\rho(x,y) + \log 4)$ may be replaced with $
c(K)\max\{\rho(x,y),\rho(x,y)^{1/K}\} \,, $ which tends to $0$ when
$x \to y\, $ and to $\rho(x,y)$ when $K \to 1\,,$ as expected.

\bigskip

\begin{lemma} {\label {anmono}} For $K>1$ the
function
$$t\mapsto\frac{2{\rm arth}(\varphi_K({\rm th}\frac{t}{2}))}{\max\{t,t^{1/K}\}}\,,$$
is monotone increasing on $(0,1)$ and decreasing on $(1,\infty)\,.$
\end{lemma}

\begin{proof} (1) Fix $K>1$ and consider
$$f(t)=\frac{2{\rm arth}(\varphi_K({\rm th}\frac{t}{2}))}{t},\quad t>0.$$  
Let $r={\rm th}\frac{t}{2}$. Then $ t/2={\rm arth}r$, and $t$
is an increasing function of $r$ for $0<r<1$. Then
$$f(t)=\frac{2{\rm arth}(\varphi_K({\rm th}\frac{t}{2}))}{t}=\frac{{\rm arth}(\varphi_K(r))}{{\rm arth}r}=F(r).$$
Then by \cite[Theorem 10.9(3)]{avvb}, $F(r)$ is strictly decreasing
from $(0,1)$ onto $(K,\infty)$. Hence $f(t)$ is strictly decreasing
from $(0,\infty)$ onto $(K,\infty)$.

(2) Next consider
$$g(t)=\frac{2{\rm arth}(\varphi_K({\rm th}\frac{t}{2}))}{t^{1/K}},$$
and let $r={\rm th}\frac{t}{2}$. Then $t=2{\rm arth}r$ and

$$g(t)=\frac{2{\rm arth}s}{2^{1/K}({\rm arth}r)^{1/K}}=
\frac{2^{1-1/K}{\rm arth}s}{({\rm arth}r)^{1/K}}\,
,$$ where
$s=\varphi_K(r)$. We next apply \cite[Theorem 1.25]{avvb}. We know
$\frac{d}{dr}({\rm arth}r)=1/(1-r^2)$.

Writing $r'=\sqrt{1-r^2}, s'=\sqrt{1-s^2},$ we obtain the quotient of the derivatives
\begin{eqnarray*}
\frac{2^{1-1/K}(1/(1-s^2))\frac{ds}{dr}}{\frac{1}{K}({\rm arth}r)^{1/K-1}(1/(1-r^2)}&=&
 2^{1-1/K}\, K \,({\rm arth}r)^{ 1-1/K}\frac{r^{'2}}{s^{'2}}\frac{1}{K}\frac{ss^{'2}\K(s)^2}{rr^{'2}\K(r)^2}\\
&=& 2^{1-1/K}({\rm arth}r)^{1-1/K}\frac{s\K(s)^2}{r\K(r)^2}
\end{eqnarray*}
by \cite[appendix E(23)]{avvb} and l'Hospital rule. By \cite[Lemma
10.7(3)]{avvb}, $\frac{\K(s)^2}{\K(r)^2}$ is increasing, since $K>1$,
$({\rm arth}r)^{1/K-1}$ is increasing. Finally, $s/r$ is increasing
by \cite[Theorem 1.25]{avvb} and E(23). So $g(t)$ is increasing in
$t$ on $(0,\infty)$.

(3) Fix $K>1$. Clearly
$$\max\{t,t^{1/K}\}=\left\{\begin{array}{lll} t^{1/K}\quad {\rm for}\quad 0\leq t\leq 1\\
t\quad {\rm for} \quad 1\leq t < \infty .\end{array}\right.$$
\noindent
Thus
$$h(t)=\frac{2{\rm arth}(\varphi_K({\rm th}\frac{t}{2}))}{\max\{t,t^{1/K}\}},\;  $$
increases on $(0,1)$ and decreases on $(1,\infty)$.
\end{proof}

\bigskip

\bigskip
\begin{figure}
\includegraphics[width=12cm]{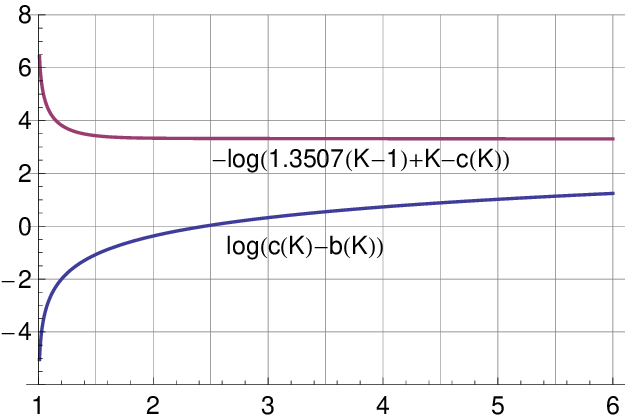}
\caption{Graphical comparison of lower and upper bounds for $c(K)$ with
 $b(K)=\log({\rm ch}(K {\rm arch}(e)))$.}
\end{figure}

\begin{subsec}{\bf Proof of Theorem \ref{refschlem}.}
{\rm  The maximum value of the function considered in Lemma \ref{anmono} is
 $c(K)=2{\rm arth}(\varphi_K({\rm th}\frac{1}{2}))$}. {\rm  The inequality
now follows from Theorem \ref{my11.2}.$\qquad\square$
}
\end{subsec}

\begin{subsec}{\bf Bounds for the constant $c(K)$.} {\rm
In order to give upper and lower bounds for $c(K)\,,$ we observe that the identity
\cite[Theorem 10.5(2)]{avvb} yields the following formula
$$c(K) = 2{\rm arth}\left( \varphi_K\left(\frac{1-1/e}{1+1/e}\right) \right)=
 2 {\rm arth} \left(\frac{1-\varphi_{1/K}(1/e)}{1+\varphi_{1/K}(1/e)} \right)\,.$$
A simplification leads to
 \begin{equation*} \label{mycK}
 c(K)=- \log \varphi_{1/K}(1/e)\,.
\end{equation*}
Next, from the inequality $ \varphi_{1/K}(r) \geq 2^{1-K}(1+r^{'})^{1-K}r^K $ for $K\ge 1, r\in (0,1)$
(cf. \cite[Corollary 8.74(2)]{avvb}) we get with $v=\log(2(1+\sqrt{1-1/e^2}))<1.3507$
\begin{eqnarray*}
 c(K)&=&- \log \varphi_{1/K}(1/e)\leq -\log \left(2^{1-K}(1+\sqrt{1-1/e^2})^{1-K}e^{-K}\right)\\
     &=& v(K-1)+K<1.3507(K-1)+K.
\end{eqnarray*}

In order to estimate the constant $c(K)$
from below we need an upper bound for $\varphi_{1/K,2}(r),\; K>1$, from above. For
this purpose we prove the following lemma.}

\begin{lemma}\label{ck} For every integer $n\geq 2$ and each $K>1,\;r\in(0,1)$, there exists
$K$-quasiconformal maps $g:\mathbf{B}^n\to \mathbf{B}^n$ and $h:\mathbf{B}^n\to\mathbf{B}^n$ with\\
$(a)\qquad g(0)=0,\;g\mathbf{B}^n=\mathbf{B}^n,\;h(0)=0,\; h\mathbf{B}^n=\mathbf{B}^n$\\
$(b)\qquad g(re_1)=\displaystyle\frac{2r^\alpha}{(1+r^{'})^\alpha+(1-r^{'})^\alpha},\;
h(re_1)=\displaystyle\frac{2r^\beta}{(1+r^{'})^\beta+(1-r^{'})^\beta}$\\
where $r^{'}=\sqrt{1-r^2}$ and $\alpha=K^{1/(1-n)}=1/\beta$.
In particular, for $n=2$ and $K>1,\; r\in(0,1)$\\
$(c)\qquad \varphi_{1/K}(r)\leq \displaystyle\frac{2r^K}{(1+r^{'})^K+(1-r^{'})^K}\;;\;\;
\varphi_{K}(r)\geq \displaystyle\frac{2r^{1/K}}{(1+r^{'})^{1/K}+(1-r^{'})^{1/K}}$.
\end{lemma}

\begin{proof} Fix $r\in(0,1)$. Let $T_a:\mathbf{B}^n\to \mathbf{B}^n$ be a
M\"obius automorphism with $T_a(a)=0$ and $T_a(\mathbf{B}^n)=\mathbf{B}^n$.
Choose $s\in(0,r)$ such that $T_{se_1}(0)=-T_{se_1}(re_1)$.
Then $\rho(0,re_1)=2\rho(0,se_1)$ \cite[(2.17)]{vuobook}, or equivalently, $(1+r)/(1-r)=((1+s)/(1-s))^2$
and hence $s=r/(1+r^{'})$. Consider the $K$-quasiconformal mapping $f:\mathbf{B}^n\to \mathbf{B}^n$,
$f(x)=|x|^{\alpha-1}x,\;\alpha=K^{1/(1-n)}$. Then $f(\pm se_1)=\pm s^\alpha e_1$. The mapping
$g=T_{-s^\alpha e_1}\circ f\circ T_{se_1}:\mathbf{B}^n\to \mathbf{B}^n$ satisfies
$g(0)=0$, $g(re_1)=te_1$ where $\rho(-s^\alpha e_1,s^\alpha e_1)=\rho(0,te_1)$ and hence
$t=2r^\alpha/((1+r^{'})^\alpha+(1-r^{'})^\alpha)$ by \cite[(2.17)]{vuobook}. The proof for $g$ is complete.
For the map $h$ the proof is similar except that we use the $K$-quasiconformal
mapping  $m:x \mapsto |x|^{\beta-1}x,\; \beta=1/\alpha$.
Note that $m=f^{-1}$ and $t=1/{\rm ch}(\alpha \;{\rm arch}(1/r))$.
For the proof of $(c)$ we apply $(a),\;(b)$ together with \cite[(3.4), p.64]{lv}.
\end{proof}

\begin{lemma} For $K>1,$ $c(K)\geq \log({\rm ch}(K {\rm arch}(e)))\geq u(K-1)+1,$ where\\
$u={\rm arch}(e){\rm th}({\rm arch}(e))>1.5412$.
\end{lemma}
\begin{proof} From Lemma \ref{ck}(c), we know that
\begin{eqnarray*}
\varphi_{1/K}(1/e)&\leq&\frac{2/e^K}{(1+\sqrt{1-1/e^2})^K+(1-\sqrt{1-1/e^2})^K}\\
                  &=&\frac{2}{(e+\sqrt{e^2-1})^K+(e-\sqrt{e^2-1})^K},
\end{eqnarray*}
hence
\begin{eqnarray*}
c(K)&=&-\log\varphi_{1/K}(1/e)\geq -\log\left(\frac{2}{(e+\sqrt{e^2-1})^K+(e-\sqrt{e^2-1})^K}\right)\\
    &=&\log\left(\frac{(e+\sqrt{e^2-1})^K+(e-\sqrt{e^2-1})^K}{2}\right)\\
    &=&\log({\rm ch}(K {\rm arch}(e)))\geq u(K-1)+1,
\end{eqnarray*}
where the last inequality follows easily from the mean value theorem, applied to
the function $p(K)= \log({\rm ch}(K {\rm arch}(e)))\,.$
\end{proof}
\end{subsec}


\end{document}